\documentclass[twoside]{amsart}

\usepackage{color}
\usepackage[colorlinks]{hyperref}
\usepackage{tikz}
\usetikzlibrary{arrows}
\usepackage{amsmath,amsthm,amsfonts,amssymb,amscd,euscript}

\newtheorem{thm}{Theorem}[section]

\theoremstyle{definition}

\theoremstyle{definition}
\newtheorem{defn}[thm]{Definition}
\theoremstyle{remark}

\numberwithin{equation}{section}

\newcommand{\filebegin}{\begin{document}}
\newcommand{\fileend}{\end{document}}
\makeatother
\def\thefootnote{}
\newcommand{\lo}{\longrightarrow}
\newcommand{\NMM}{\hspace*{2mm}}

\renewcommand{\baselinestretch}{1.1}

\renewcommand{\baselinestretch}{1.1}
\def\n{\noindent}%

\numberwithin{equation}{section}
\def\mapdown#1{\Big\downarrow\rlap
{$\vcenter{\hbox{$\scriptstyle#1$}}$}}
\begin{document}

\title{On half-synchronized systems}
\author{ S. Jangjooye Shaldehi }%
\address{Faculty of Mathematical Sciences, Alzahra University }
\email{ sjangjoo90@gmail.com}%


\begin{abstract}
 
 A subclass of coded systems containing synchronized systems is the family of half-synchronized systems. In this note, we will consider them and show that the property  `half-synchronized' lifts under hyperbolic maps. This enables us to define an equivalence relation on the set of
half-synchronized systems. Also,  when the domain is half-synchronized, we show that right-closing a.e., 1-1 a.e. factor maps have a strong decoding property.

\end{abstract}

\maketitle

\section{Introduction}

Coded systems were defined by Blanchard and Hansel \cite{BH} as a generalization of sofic shifts. Amongst subshifts,
a well-known subclass of coded systems is the family of synchronized systems. The purpose of this paper is to generalize parts of the synchronized theory to certain subclass of coded systems, the so-called half-synchronized systems. The half-synchronized systems  often serve as a landmark within the coded systems, showing the difficulties in extending the class of synchronized systems and keeping a satisfactory analogon of the sofic shifts theory.

 The study of conjugacy  is of interest in symbolic dynamics. But there is no general algorithm for deciding whether two
shift spaces are conjugate. Thus it makes sense to ask if there is an equivalence relation, weaker than conjugacy, on some subclasses of subshifts. The investigation for the existance of common extensions where
the factor maps have certain properties in common is of interest in coded systems and  has a long history. Amongst the properties, we are interested in hyperbolicity. Common extensions with hyperbolic maps lead to a version of almost conjugacy for coded systems. Fiebig in \cite{F2} shows that   having
a common coded (or synchronized) hyperbolic  extension is an equivalence relation on the set of
coded (or synchronized) systems. In Section \ref{two}, first we prove that  `being half-synchronized' is an invariant for hyperbolic maps (Theorem \ref{half}). Then, We show that a common half-synchronized hyperbolic   extension is an equivalence relation on the set of
half-synchronized systems (Theorem \ref{equivalency}).

Closing maps are important in the
coding theory \cite{Kit, LM}. They also have a natural
description in hyperbolic dynamics
 \cite{BS}: right-closing (resp.  left-closing) maps are injective on
unstable (resp. stable)  sets. When the domain is a synchronized system, right-closing a.e., 1-1 a.e. factor maps have a strong decoding property \cite{F3}. This property gives finitary regular isomorphism between general shift spaces. In section \ref{three}, we generalize this result to half-synchronized systems (Theorem \ref{decoder}).

\section{Background and Notations}

Let $\mathcal A$ be a non-empty finite set. The \emph{full $\mathcal
A$-shift}
 denoted by ${\mathcal A}^{\mathbb Z}$, is the set of all bi-infinite sequences of symbols from ${\mathcal A}$.
A \emph{block} over ${\mathcal A}$ is a finite sequence of
symbols from ${\mathcal A}$. The \emph{shift map} on ${\mathcal A}^{\mathbb Z}$ is the map $\sigma$ where
$\sigma(\{x_i\})=\{y_i\}$ is defined by $y_{i}=x_{i+1}$. The pair $(\mathcal {A}^{\mathbb Z},\,\sigma)$ is the \emph{full shift} and any closed invariant subset of that is called a \emph{shift space}.

Denote by ${\mathcal B}_{n}(X)$  the set of all admissible $n$-blocks and let ${\mathcal B}(X)=\bigcup_{n=0}^{\infty}{\mathcal B}_{n}(X)$
be  the \emph{language} of $X$.
For $ u\in\mathcal B(X) $, let the \emph{cylinder}
$ [u] $ be the set $ \{x\in X:\,x_{[l,l+|u|-1]}=u\} $. 

Suppose $X$ is a
subshift over alphabet ${\mathcal A}$. For $m,\,n\in \mathbb Z$ with $-m
\leq n$, define the \emph{$(m+n+1)$-block map} $\varphi: {\mathcal
B}_{m+n+1}(X) \rightarrow {\mathcal D}$ by
\begin{equation}\label{2.1} 
y_{i}=\varphi(x_{i-m}x_{i-m+1}...x_{i+n})=\varphi(x_{[i-m,i+n]})
\end{equation}
where $y_{i}$ is a symbol in alphabet ${\mathcal D}$. The map $\varphi=\varphi_{\infty}^{[-m,n]}: X
\rightarrow {\mathcal D}^{\mathbb Z}$ defined by $y=\varphi(x)$ with
$y_{i}$ given by \ref{2.1} is called the \emph{sliding block code} (or \emph{code}) induced by
$\varphi$. So there is $ N\geq 0 $ such that $ \varphi(x)_0 $ is determined by $ x_{[-N,N]} $. We call $ 2N+1 $ a coding length for $ \varphi $. If $ m=n=0 $, then $ \varphi $ is \emph{$ 1 $-block code} and $ \varphi=\varphi_{\infty} $. An onto (resp. invertible) code is called a
\emph{factor map} (resp. \emph{conjugacy}).

A point $  x\in  X$ is \emph{doubly transitive} if every
block in $  X$ appears in $  x$ infinitely often to the left and to the right. Let $ \varphi:X\rightarrow Y $ be a factor map. If
there is a positive integer $  d$ such that every doubly transitive point of $  Y$
has $  d$ pre-images, then we call $  d$ the \emph{degree} of $ \varphi $.

 Let $G$ be a directed graph and $  \mathcal V$ (resp. $\mathcal E$) the set of its vertices (resp. edges) which is supposed to be countable.
 An \emph{edge shift}, denoted by $X_{G}$,  is a shift space which consist of all
  bi-infinite sequences of edges from $ \mathcal E $.

A labeled graph ${\mathcal G}$ is a pair $(G,{\mathcal L})$ where
$G$ is a graph and  
${\mathcal L}: {\mathcal E} \rightarrow {\mathcal A}$ its labeling. Associated to $\mathcal G$, a space 
$$ X_{{\mathcal G}}=\text{closure}\{{\mathcal L}_{\infty}(\xi): \xi \in X_{G}\}=\overline{{\mathcal L}_{\infty}(X_{G})}$$
is defined and ${\mathcal G}$ is called a  \emph{cover} of $X_{{\mathcal
G}}$. When $G$ is a finite graph, $X_{{\mathcal G}}={\mathcal L}_{\infty}(X_{G})$ is a \emph{sofic shift}.

A block $v \in {\mathcal B}(X)$ is called \emph{synchronizing} if whenever $uv,\,vw\in{\mathcal B}(X)$,
we have $uvw \in {\mathcal B}(X)$. An
irreducible subshift  $X$ is a \emph{synchronized system} if it has a synchronizing block. 

A \emph{coded} system is the closure of the set of sequences obtained by freely concatenating
the blocks in a list of blocks.

A labeled graph ${\mathcal G}=(G,{\mathcal L})$ is called \emph{right-resolving} if for each vertex $I$ of $G$, the edges starting at $I$ carry
different labels. A \emph{minimal right-resolving cover} of a sofic shift $X$ is a right-resolving cover
 having the fewest vertices among all right-resolving covers of $X$. It is unique up to isomorphism \cite[Theorem 3.3.18]{LM} and is called the \emph{Fischer cover} of $ X $.

For $ x\in\mathcal B(X) $, call  $ x_{-}=(x_{i})_{i<0} $ (resp. $ x_{+}=(x_{i})_{i\in\mathbb Z^{+}} $) the \emph{left (resp. right) infinite $ X $-ray} and let $ X^{+}=\{x_{+}:\,x\in X\} $. The follower set of  $ x_- $ is defined as $ \omega_{+}(x_{-})=\{x_{+}\in X^{+}:\,x_{-}x_{+}$ is a point in $X\} $.

\section{One equivalence relation for half-synchronized systems}\label{two}

We prove that hyperbolic maps lift the property half-synchronized. We use this
to show that common half-synchronized hyperbolic extensions define an equivalence
relation on the set of half-synchronized systems.

\begin{defn}\label{defnhalf} 
 A transitive subshift $  X$ is \emph{half-synchronized}, if there is a
block $  m\in\mathcal B(X)$ and a left-transitive point $  x\in X$ such that $ x_{[-|m|+1,0]}=m $ and $ \omega_{+}(x_{(-\infty,0]})=\omega_{+}(m) $. Then, $  m$ is called a \emph{half-synchronizing block} for $  X$.
\end{defn}  
Dyck shift and synchronized systems are half-synchronized. Now we review the concept of the Fischer cover for a half-synchronized system \cite{FF}.
Let the collection of all follower sets $ \omega_{+}(x_{-})$ be the set of vertices of a graph $ X^{+} $. There is an edge from $ I_{1} $ to $ I_{2} $ labeled $ a $ if and only if there is a $ X $-ray $ x_{-} $ such that $ x_{-}a $ is a $ X $-ray and $I_{1}=\omega_{+}(x_{-}) $, $ I_{2}=\omega_{+}(x_{-}a) $. This labeled graph is called the \emph{Krieger cover} for $ X $.
If $ X $ is a half-synchronized system with half-synchronizing block $ \alpha $, the irreducible component of the Krieger cover containing the vertex $ \omega_{+}(\alpha)  $ is called the Fischer cover of $ X $.

The following definition is motivated by the notion of hyperbolic
homeomorphism between the sets of doubly transitive points introduced in \cite{T}.
\begin{defn}\label{hyper} 
 Let $  X$ and $  Y$ be transitive subshifts. The factor map $ \varphi:X\rightarrow Y $ is
\emph{hyperbolic} if there is a $ d\in \mathbb N $ and a block $  w$ and $  d$ blocks
$$ m^{(1)},m^{(2)},\ldots,m^{(d)}\in\mathcal B_{2k+1}(X), $$
such that
\begin{enumerate}
\item
 if $ y\in Y $ such that $ y_{[-n,n]}=w $, then $$ \varphi^{-1}(y)_{[-k,k]}=\{m^{(1)},m^{(2)},\cdots,m^{(d)}\}. $$
\item
if $w'\in\mathcal B(Y)  $ beginning and ending with $  w$, then for each $ 1\leq i\leq d $ there is a
unique block $ a^i\in\mathcal B(X) $, such that for any $ x\in X $ with $  \varphi(x)_{[-n,n+p]}=w'$ and
$ x_{[-k,k]}=m^(i) $, it holds $ x_{[-k,k+p]}=a^i $.
\end{enumerate}
\end{defn}

The next theorem say that \cite[Theorem 4.2]{F2} that
were stated for synchronized systems are actually valid for half-synchronized systems as well. In fact, it says that the property  `half-synchronized' is weak enough to
lift under hyperbolic maps. 
\begin{thm}\label{half} 
Let $  X$ and $  Y$ be transitive subshifts and $\varphi:X\rightarrow Y  $ a hyperbolic factor
map. Then, $  X$ is half-synchronized if and only if $  Y$ is half-synchronized.
\end{thm}
\begin{proof}
First suppose $ X  $ is half-synchronized and $ \mathcal G=(G,\mathcal L) $ is the Fischer cover of $ X $. By \cite[Theorem 1.4]{FF}, $ \mathcal G $ has residual image in the one-sided shift. Since $\varphi=\Phi_{\infty}  $ has a degree \cite[Theorem 3.2]{F2}, $ (G,\Phi\circ\mathcal L) $ has  residual image in $ $  $Y^+ $. So $ Y $ is half-synchronized.

Now let $  Y$ be half-synchronized. Since $ \varphi $ is hyperbolic, it has a degree $  \ell$ \cite[Theorem 3.2]{F2}.  Hyperbolicity of $ \varphi $ implies that the transitive left ray $ y_- $ does not have less than $  \ell$ preimages. On the other hand, if it has more than $  \ell$ preimages, then a prolongation of $ y_- $ gives a point in $ Y $ with more than $  \ell$ preimages. So $ y_- $ has $  \ell$  preimages.

Let $  w\in\mathcal B(Y)$ and $ m^{(1)},m^{(2)},\cdots,m^{(d)}\in\mathcal B(X) $  satisfying Definition \ref{hyper}. Without loss of generality, we may assume that $ \varphi $ is $ 1 $-block and $  w$ is half-synchronizing block. So there is a transitive  left ray $y_-$   terminating at $w$ such that $ \omega_{+}(w)=\omega_{+}(y_-) $. Let $  m\in\mathcal B_{2n+1}(X)$ such that $ \Phi(m)=w $ and $ x_{1_-},x_{2_-},\cdots ,x_{\ell_-}$ be $  \ell$ preimages for $ y_- $.  Then, since $ \varphi $ has a degree, $ x_{1_-} $ is transitive \cite[Theorem 3.2]{F2}. We may assume that  $ m=m^{(1)} $ and $ x_{1_-} $ terminates at $m$. Now we show that $ m $ is a half-synchronizing block for $ X $. 

Suppose $ mu\in\mathcal B(X) $ and prolongate $ mu $ such that $ muvm\in\mathcal B(X)$. Then, $ w\Phi(uv)w\in\mathcal B(Y) $ and since $ w$ is half-synchronizing, $\Phi(uv)w\in \omega_{+}(y_-) $. Now $ y_{-}\Phi(uv)w $ also has $ \ell $ preimages $ x_{i_-} u_{i}$ for $ 1\leq i\leq \ell $. (2) of Definition \ref{hyper} implies that $ u_1=uv$. So $ u\in \omega_{+}(x_{1_-}) $.

\end{proof}
Theorem \ref{half} enables  us to introduce an equivalence relation
using hyperbolic factor codes on the set of half-synchronized systems.
\begin{thm}\label{equivalency} 
 Having a common hyperbolic extension is an equivalence relation on
the half-synchronized systems.
\end{thm}
\begin{proof}
Let  $ X $  and $ Y $ (resp.  $ Y $  and $ Z $) have a common half-synchronized hyperbolic extension $ (V,\varphi_{X},\varphi_{Y})$ (resp. $ (W,\varphi'_{Y},\varphi'_{Z}) $). Suppose $ (\Sigma,\psi_{V},\psi'_{W})$ is the fiber product of $ (\varphi_{Y},\varphi'_{Y}) $ and $  \Gamma  $ is an irreducible component of $ \Sigma $ such that the restriction of $ \psi_{V} $ and $ \psi'_{W} $ to $ \Gamma $  are onto. Since $ \varphi_{Y} $ is hyperbolic, $ \psi'_{W}:\Gamma\rightarrow W $ is hyperbolic and so by Theorem \ref{half}, $\Gamma$ is half-synchronized. The result then
follows from the fact that compositions of
hyperbolic maps are hyperbolic.
\end{proof}

\section{The existence of decoder block}\label{three}

A factor code $ \varphi:X\rightarrow Y $ is 1-1 a.e. if any doubly transitive point in $  Y$ has exactly one preimage. 
 It is  \emph{right-closing almost everywhere}  if there is $ n\in\mathbb N $ such that for any two left-transitive points $ x,y\in X $ 
with $ x_{(-\infty,0]}=y_{(-\infty,0]} $ and $ \varphi(x)_{(-\infty,n]}=\varphi(y)_{(-\infty,n]} $, it holds $ x_{1}=y_{1} $. Then, $ \varphi $ is called \emph{$  n$-step right-closing a.e.}. 
  
\begin{defn}
Let $ \varphi:X\rightarrow Y $ be a factor map between arbitrary subshifts $ X $ and $ Y $. A block $ w\in\mathcal B(Y) $ is a \emph{decoder block} for $ \varphi $ if there is $ k\in\mathbb N $, the \emph{anticipation} of $ w $, such that for all $ n\in\mathbb N $ and all points $ x,y\in X $ with $ \varphi(x)_{[-|w|+1,0]}=\varphi(y)_{[-|w|+1,0]}  $ and $ \varphi(x)_{[1,n+k]}=\varphi(y)_{[1,n+k]}  $, it holds that $ x_{[1,n]}=y_{[1,n]} $.
\end{defn}
Remind that if $ \varphi:X\rightarrow Y $ is a factor map with a decoder block and $ \mu $ an ergodic measure on $ X $ with full support, then  $\varphi:(X,\mu)\rightarrow (Y,\nu)  $ is a finitary regular isomorphism, where $ \nu=\mu\circ\varphi^{-1} $ \cite{F3}.

Fiebig shows that when the domain is a synchronized system, right-closing a.e., 1-1 a.e. factor maps have a strong decoding property \cite{F3}. We generalize this property to half-synchronized systems.
The ingredients for the proof is almost similar to \cite[Theorem 2.3]{F3}.
\begin{thm}\label{decoder} 
Let $ X $ be a half-synchronized system and  $\varphi:X\rightarrow Y  $ a factor map. Then, $ \varphi $ is right-closing a.e., 1-1 a.e. if and only if it has a decoder block.

\end{thm}
\begin{proof}
Suppose $ \varphi $ is $ k $-step right-closing a.e. and $  m\in\mathcal B(X)  $ is a half-synchronizing block with length greater than a coding length for $ \varphi $. Let $ z\in X $ be a doubly transitive point satisfying Definition \ref{defnhalf}. By compactness combined with the fact that $ \varphi $ is 1-1 a.e., there is $ p\in \mathbb N $ such that for
 $ x\in X $ with $ \varphi (x)_{[-p,p]}=\varphi(z)_{[-p,p]} $ it holds that $ x_{[-|m|+1,0]}=m $. We claim that $ w:= \varphi(z)_{[-p,p]}$ is a decoder block with anticipation $ k $.

Let $ x,y\in X $ with $ \varphi(x)_{[-|w|+1,0]}=\varphi(y)_{[-|w|+1,0]}=w  $ and $ \varphi(x)_{[1,n+k]}=\varphi(y)_{[1,n+k]}  $. Then, $ x_{[-p-|m|+1,-p]}=y_{[-p-|m|+1,-p]}=m $. Since $ m $ is half-synchronizing, there are points $ x',y'\in X $ as,
$$\hspace{-40mm}x'_{i}=\left\{
\begin{tabular}{ll}
$z_{i+p}$&$i\leq -p,$ \\
$x_{i}$&$i\geq -p+1,$ \\
\end{tabular}
\right .$$
\vspace{-12mm}
$$\hspace{60mm}y'_{i}=\left\{
\begin{tabular}{ll}
$z_{i+p}$&$i\leq -p,$ \\
$y_{i}$&$i\geq -p+1.$ \\
\end{tabular}
\right .$$
$ x' $ and $ y' $ are left transitive and $ \varphi(x')_{(-\infty,n+k]}=\varphi(y')_{(-\infty,n+k]} $. Since $ \varphi $ is right-closing a.e., $  x'_{[-p+1,n]}=y'_{[-p+1,n]}$. So $ x_{[1,n]}=y_{[1,n]} $.

 Now suppose $  w$ is a decoder block with anticipation $  k$ and $ y\in D(Y) $. So $ w $ appears in
$  y$ infinitely often to the left and to the right. Since $  w$ is a decoder block, $ \varphi$ is 1-1 
a.e.. If $ y_{[i-|w|+1,i]}=w $ for some $ i\in\mathbb Z $, a coordinate $ j>i $ of the preimages of $  y$ is determined by $ y_{[i-|w|+1,j+k]} $. Therefore, $ \varphi $ is $ k $-step right-closing a.e.. □ 
\end{proof}


\providecommand{\bysame}{\leavevmode\hbox
to3em{\hrulefill}\thinspace}


\end{document}